\documentclass[11pt,amscd]{amsart}
\usepackage{amsaddr}
\usepackage{graphicx,amssymb,amsmath,amsthm}

\usepackage[utf8]{inputenc}
\usepackage{epstopdf}
\usepackage[boxed,linesnumbered]{algorithm2e}	
\usepackage{bm}
\usepackage[numbers]{natbib}
\usepackage{color}
\usepackage{listings}
\usepackage{array}

\usepackage[margin=1.25in]{geometry}
\newcolumntype{L}{>{$}l<{$}} 
\newcolumntype{C}{>{$}c<{$}} 

\usepackage{bm,graphicx,textcomp,bbm} 

\usepackage{booktabs}

\usepackage{subcaption}
\usepackage{lipsum}
\usepackage{amsfonts}
\usepackage{graphicx}
\usepackage{epstopdf}
\usepackage{algorithmic}
\ifpdf
  \DeclareGraphicsExtensions{.eps,.pdf,.png,.jpg}
\else
  \DeclareGraphicsExtensions{.eps}
\fi
\newcommand{\eps}{{\varepsilon}}

\newcommand{\R}{\mathbb{R}}

\def\O{{\mathcal{O}}}
\newcommand{\dotx}{{\dot x}}
\newcommand{\doty}{{\dot y}}
\newcommand{\x}{x}
\newcommand{\y}{y}
\newcommand{\E}{\mathbb{E}}
\newcommand{\V}{\mathbb{V}}
\newcommand{\C}{\mathcal{C}}
\newcommand{\eff}{v}
\newcommand{\kapfactor}{\kappa}

\usepackage{amsopn}

\allowdisplaybreaks

\title{A note on statistical consistency of numerical integrators for multi-scale dynamics} 
\author{J. Frank$^\star$ and G.~A.~Gottwald$^\dagger$}
\address{$^\star$Mathematical Institute, Utrecht University, P.O. Box 80010, 3508 TA Utrecht, the Netherlands}
\address{$^\dagger$ School of Mathematics and Statistics, University of Sydney, NSW 2006, Australia}
\email[J. Frank and G. A. Gottwald]{J.E.Frank@uu.nl {\rmfamily and} georg.gottwald@sydney.edu.au}

\begin{document}

\maketitle 

\begin{abstract}
A minimal requirement for simulating multi-scale systems is to reproduce the statistical behavior of the slow variables. In particular, a good numerical method should accurately aproximate the probability density function of the continuous-time slow variables. In this note we use results from homogenization and from backward error analysis to quantify how errors of time integrators affect the mean behavior of trajectories. We show that numerical simulations converge, not to the exact probability density function (pdf) of the homogenized multi-scale system, but rather to that of the homogenized modified equations following from backward error analysis. Using homogenization theory we find that the observed statistical bias is exacerbated for multi-scale systems driven by fast chaotic dynamics that decorrelate insufficiently rapidly. This suggests that to resolve the statistical behavior of trajectories in certain multi-scale systems solvers of sufficiently high order are necessary. Alternatively, backward error analysis suggests the form of an amended vector field that corrects the lowest order bias in Euler's method.  The resulting scheme, a second order Taylor method, avoids any statistical drift bias. We corroborate our analysis with a numerical example.
\end{abstract}

\maketitle


\section{Introduction}
When simulating complex multi-scale dynamics one is often interested in the accurate description, not of the full system with all its degrees of freedom, but of only some distinct relevant variables, for example the slow variables. Numerical weather forecasting provides a good example, where we are interested in the dynamics of the large-scale high and low pressure fields which evolve on time scales of days rather than in fast buoyancy oscillations of the atmosphere's stratification surfaces. Another example is decadal climate prediction where we are not interested in the actual dynamics of the large-scale atmospheric weather but rather in their effect on the slowly evolving oceanic patterns such as El-Ni\~no. Whereas in the example of numerical weather prediction we desire accurate time evolution of the slow relevant variables, in climate science we are often more interested in statistical properties such as mean global temperature or the frequency of extreme events. Reproducing such mean statistical behavior of the slow variables is a minimal requirement for any simulation of complex multi-scale systems. It has long been recognized that the numerical discretization scheme employed to simulate a dynamical system profoundly affects the numerically observed statistical behavior \cite{AbKoMa03,DubinkinaFrank07,DubinkinaFrank10}. This note is concerned with the problem of numerically integrating multi-scale systems with the aim to reliably recover their statistical properties.\\ 

We consider here deterministic multi-scale systems of the form
\begin{align}
 \dotx & = \frac{1}{\varepsilon}\, h(\x)f_0(\y)+ f(\x,\y), \quad \x(0)=\xi \label{e.ms.1}\\
\doty & = \frac{1}{\varepsilon^2}\, g(\y), \quad \y(0)=\eta, \label{e.ms.2}
\end{align}
with $\x\in\R^d$, $\y\in\R^\ell$. 
The parameter $\eps\ll 1$ characterizes the degree of time scale separation. Here the slow dynamics evolves on a characteristic time of order $1$ and the fast dynamics on a characteristic time of $\varepsilon^2$.  
We assume that the vector fields $f_0:\R^\ell\to\R^m$, $h:\R^d\to L(\R^d,\R^m)$, $f:\R^d\times\R^\ell\to\R^d$ and $g:\R^\ell\to\R^\ell$ satisfy certain regularity conditions and that the fast $y$-dynamics is sufficiently chaotic with compact chaotic attractor $\Lambda\subset\R^\ell$ and ergodic invariant probability measure $\mu$. We consider the case when $\int_\Lambda f_0\,d\mu=0$, i.e.~when classical averaging would yield trivial constant-in-time dynamics. In this situation the slow dynamics exhibits stochastic dynamics on the slow time scale ${\mathcal{O}}(1)$ \cite{Givonetal04,PavliotisStuart}.\\ 

Numerical simulation of the multi-scale system (\ref{e.ms.1})--(\ref{e.ms.2}) is challenging: To capture the slow dynamics of interest, for any fixed value of the time scale separation parameter $\varepsilon$, we obtain convergence in the limit $\Delta t \to 0$, but for $\varepsilon$ small, the time step $\Delta t$ used to propagate the slow variables must be chosen of the order of $\varepsilon^2$ to resolve the fast dynamics and meet stability restrictions, making direct numerical simulations computationally impractical. A minimal requirement for a numerical integrator is that it should reproduce the statistical behavior of the slow variables of interest. Ideally we would like to employ $\Delta t\sim\mathcal{O}(1)$\footnote{ 
Special multi-scale methods have been devised to do this (see, e.g.~\cite{E03,EEtAl07,GearKevrekidis03,KevrekidisGearEtAl03}).}. 
However, we will see that depending on the statistical behavior of the fast dynamics, in particular on the decay of the correlation function of $f_0(y)$, a time step $\Delta t \sim \mathcal{O}(\varepsilon^2)$ may not be sufficient to recover even the statistical behavior of the slow dynamics and one will need time steps such that $\kappa=\Delta t/\varepsilon^2 \to 0$ as $\eps\to 0$. (Note that for $\kappa>0$ solutions of the fast integrator do not converge in the limit $\varepsilon\to 0$ to the exact solution of \eqref{e.ms.2}.)  
In other words, it is insufficient to simply resolve the fast motions as $\varepsilon\to 0$, one must in fact \emph{accurately approximate them} in this limit, even when the goal is to determine the mean behavior of the slow variables.
This inability of numerical time steppers of order $p$ to reproduce the statistical behavior of the slow dynamics will be linked to the persistence of $\O(\Delta t^p)$-terms in the backward error analysis; furthermore these error terms have a quantifiable influence on the long-time statistics as they will be shown to correspond to drift corrections in the homogenized diffusive limit equations of the numerical discrete time maps.\\

For multi-scale systems of the form (\ref{e.ms.1})--(\ref{e.ms.2}) the statistical behavior of the slow dynamics, in the limit of infinite time scale separation $\eps\to 0$, is described by a stochastic differential equation (SDE) which can be explicitly stated. The mathematical tool to describe the long-time stochastic behavior of slow dynamics is known as {\em{homogenization}}  \cite{Givonetal04,PavliotisStuart}. Homogenization describes the integrated effect of the fast (either stochastic or chaotic) dynamics on the slow variables as noise. Initially developed for stochastic multi-scale systems \cite{Khasminsky66,Kurtz73,Papanicolaou76}, homogenization has been extended recently to deterministic multi-scale systems. In the deterministic case the theory is restricted to the skew-product case (\ref{e.ms.1})--(\ref{e.ms.2}) in which the slow dynamics does not couple back to the fast dynamics. The fully coupled case poses the potential problem that the invariant measure of the fast dynamics may not vary smoothly with the slow variable; in this instance the averaged vector fields may not even be Lipshitz and uniqueness and existence of the homogenized equation may not be guaranteed. For the deterministic skew product case (\ref{e.ms.1})--(\ref{e.ms.2}), it was shown rigorously that for sufficiently chaotic fast dynamics the emergent stochastic long-time behavior of the slow dynamics is given by stochastic differential equations driven by Brownian motion \cite{MelbourneStuart11,GottwaldMelbourne13c,KellyMelbourne17}. The assumed mild conditions on the chaoticity of the fast $y$-dynamics are satisfied by a large class of maps and flows. 
For maps, the convergence to Brownian motion holds when the correlation function is summable. For flows, it suffices that there is a Poincar\'e map with these properties (irrespective of the mixing properties of the flow). These include, but go far beyond, Axiom A diffeomorphisms and flows, H\'enon-like attractors and Lorenz attractors.  Precise statements about the validity can be found in \cite{MelbourneNicol05,MelbourneNicol08,MelbourneNicol09}. We remark that for weakly chaotic dynamics when the correlations are not summable, the noise is not Brownian anymore but rather $\alpha$-stable \cite{GottwaldMelbourne13c}\footnote{We use the terminology strongly and weakly chaotic here in a manner different from the usual distinction between exponential and algebraic decay of correlations; cf.\cite{GottwaldMelbourne13}.}. Homogenization has been used as a framework for stochastic parametrizations in the context of numerical weather forecasting and climate science \cite{MTV99,Majdaetal01,MTV02,Majda03,monahan_stochastic_2011,culina_stochastic_2011,GottwaldEtAl17} and is at the core of the design of several efficient numerical multi-scale integrators such as the {heterogeneous multi-scale method} \cite{E03,EEtAl07} and {equation-free} projection \cite{GearKevrekidis03,KevrekidisGearEtAl03}.\\

Depending on the underlying deterministic dynamical multi-scale system, the noise appearing in the limiting homogenized SDE can be either additive or multiplicative. It is well known that the solution of an SDE is sensitive to the approximation of the Brownian motion. This sensitivity gives rise to the different interpretations of the noise such as It\^o versus Stratonovich interpretations (see the insightful discussion in \cite{HorsthemkeBook}).
In \cite{GottwaldMelbourne13c} it was shown that in the case when the slow dynamics is one-dimensional the stochastic differential equation describing the diffusive behavior of the slow dynamics is to be interpreted in the Stratonovich sense. The intuitive argument for this result is that the noisy SDE is a rough approximation of a smooth dynamical system, hence in the limiting process of infinite time scale separation classical calculus should prevail which necessitates the Stratonovich interpretation\footnote{This does not hold for higher-dimensional slow sub-spaces where the noise is neither Stratonovich nor It\^o \cite{KellyMelbourne17} and the conditions for the Wong-Zakai theorem are not satisfied.}. The limiting SDE for deterministic discrete-time maps, however, was shown to be neither of Stratonovich nor of It\^o type. The noise is It\^o only if the fast dynamics is $\delta$-correlated.\\ 

This immediately points to a problem when numerically simulating a continuous-time multi-scale system: The long-term statistics of a dynamical multi-scale system, be it continuous-time or discrete time, is described by its homogenized limiting SDE. However, the limiting stochastic differential equation describing the long-time statistical behavior of the discretized slow dynamics, that is of the numerical integrator, might be different from that of the continuous-time system it is designed to model. Using backward error analysis, we show that the leading-order term responsible for the difference is the limiting second order contribution of the modified equation corresponding to the numerical map. The main contribution of our work is to show that the local errors of a time stepper generate a long-time error of the mean behavior which is recovered by homogenization theory. These error terms are of the order $\O(\Delta t^p)$ for a $p$th order integrator. This result allows us to draw an important practical conclusion: In order for a numerical discretization scheme to reproduce the long-time statistical behavior of the slow dynamics it may be necessary to employ a sufficiently high order time-stepping method. In particular, the Euler scheme can lead to massively different statistical behavior with strong bias. This is the case when, as we will see, the fast chaotic dynamics does not decay sufficiently quickly and its statistical behavior is far from being close to independent identically distributed ({\em{i.i.d.}}) random variables. In contrast, first order schemes \emph{are} sufficient to capture the long-time statistical behavior for multi-scale systems with chaotic fast dynamics exhibiting rapid decay of correlation. 
As we will see, discretization-induced biases can be expressed using homogenization theory. This allows us to explicitly subtract the bias from the slow vector field of the deterministic equation (\ref{e.ms.1}), resulting in a remarkably accurate explicit time stepper.\\

The paper is organized as follows. In Section~\ref{s.ms} we introduce the diffusive limit of the deterministic multi-scale system (\ref{e.ms.1})--(\ref{e.ms.2}) and of its associated Euler scheme. The diffusive limits of the original continuous-time deterministic multi-scale system and its Euler discretized version are shown to differ in the drift term. In Section~\ref{s.bea} we present the backward error analysis of Euler's method and Heun's method and the homogenized limit of the lowest order modified equation for each, describing how its respective long-time statistics differs from that of (\ref{e.ms.1})--(\ref{e.ms.2}). Section~\ref{s.numerics} presents numerical simulations corroborating our analytical results. We conclude with a summary and an outlook in Section~\ref{s.summary}.

\section{The diffusive limit of the multi-scale system and its Euler scheme}
\label{s.ms} 
Using fairly weak conditions on the chaoticity of the fast $y$ dynamics, it was recently proved  in \cite{MelbourneStuart11,GottwaldMelbourne13c,KellyMelbourne17} that the long-term behavior of deterministic multi-scale systems (\ref{e.ms.1})--(\ref{e.ms.2}) is stochastic and is described on times of order ${\mathcal{O}}(1)$ by the following homogenized stochastic differential equation
\begin{align} 
\label{e.SDE_flow}
dX=F(X)\,dt + \sigma h(X)\circ dW_t, \quad X(0)=\xi.
\end{align}
For simplicity of exposition and ease of computation, we choose in the following $d=m=1$. The drift term is given by $F(X)=\int_\Lambda f(X,y)\,d\mu$, $W_t$ is unit $1$-dimensional Brownian motion with the variance given by a Green-Kubo formula with
\begin{align}
\frac{1}{2}\sigma^2 = \int_0^\infty \C[f_0(y)](t) \, dt\, ,
\label{e.GK}
\end{align}
where $\C[f_0(y)](t) = \E [f_0(y)f_0(\varphi^ty)]$ denotes the autocorrelation function of $f_0$ with $\varphi^t$ denoting the flow of the vector field $g(y)$ (in particular, $\varphi^t$ is independent of $\varepsilon$), and the expectation 
\[
	\E [A] = \int_\Lambda A(y) d\mu
\]
is taken with respect to the fast invariant measure $\mu$. As discussed in the Introduction, the noise is of Stratonovich type because the smooth dynamical system \eqref{e.ms.1}--\eqref{e.ms.2} is approximated by a rough SDE \eqref{e.SDE_flow}, and hence classical calculus has to be valid throughout the limiting procedure of homogenization. For the precise statements we refer the interested reader to \cite{GottwaldMelbourne13c}.\\

When the multi-scale system \eqref{e.ms.1} is discretized with time step $\Delta t$  by a numerical integration method, the slow dynamics is given by a map.  For instance, the first order forward Euler method gives
\begin{align}
x_{n+1}=x_n+\Delta t\frac{1}{\varepsilon}\, h(x_n)f_0(y_n)+\Delta t \,f(x_n,y_n), \quad x_0=\xi,
\label{e.map}
\end{align}
where $y_n \approx y(n\Delta t)$ is also obtained via a map $y_{n+1} = \Phi(y_n)$ that approximates $\varphi^t$ on time $\Delta t$. In this paper, we compute $\Phi(y)$ using multiple time stepping, through the $K$-fold application of the same numerical integrator as used for the slow dynamics\footnote{An alternative strategy would be to use the method suggested in \cite{TaoEtAl10,ArielEtAl13}.},
\begin{equation}\label{disc.fast}
	y_{n,k+1} = y_{n,k} + \delta t\, \varepsilon^{-2} g(y_{n,k}), \quad k=0, \dots, K-1,
\end{equation}
with initial condition $y_{n,0} = y_n$
and time step $\delta t = \Delta t/K$. 
We set $y_{n+1} =  y_{n,K}$ to define the map $y_{n+1} = \Phi(y_n)$. 
In the limit $\varepsilon\to 0$ we choose the scaling $\Delta t =\kapfactor\,\varepsilon^2$, where $\kapfactor>0$ is a small but finite constant (i.e. we solve the slow equation on the fast time scale). This implies that the effective stepsize of the fast motion in \eqref{disc.fast} is $\delta t \,\varepsilon^{-2} = \kapfactor/K$ and the map $\Phi$ is independent of $\varepsilon$.  Consequently, the fast motion \eqref{disc.fast} does not converge in the limit $\varepsilon\to 0$ to the exact solution of \eqref{e.ms.2}.  Instead, the constant $K$ is chosen such that the fast motion is well-resolved for all $\varepsilon$.  We also assume that the discrete dynamics \eqref{disc.fast} possesses a chaotic attractor that satisfies the conditions needed for the existence of the SDE limit as discussed below.

For the map (\ref{e.map}) it was rigorously proven in \cite{GottwaldMelbourne13c} that the long-time statistics on times of order ${\mathcal{O}}(1/\eps^2)$ is governed by the following SDE
\begin{align} 
\label{e.SDE_map0}
dX=\Bigl(\kappa\,F(X)-\frac12 \kappa^2 h(X)h'(X)\,\E[ f_0^2] \Bigr)\, dt + \kappa\, \hat\sigma h(X)\circ dW_t\, ,
\quad X(0)=\xi,
\end{align}
where $F(X)$ is the same as before for the continuous-time system, $W_t$ is again unit $1$-dimensional Brownian motion and the variance is given by a Green-Kubo formula
\begin{align} 
\label{e.GK_disc}
\hat\sigma^2=\E[ f_0^2]+
\sum_{n=1}^\infty \E[ f_0(\Phi^n\,y)f_0(y)] +
\sum_{n=1}^\infty \E[f_0(y)f_0(\Phi^n\,y)],
\end{align}
where $\Phi^n$ denotes the $n$-fold application of the discrete map $\Phi$.
Evaluating the time integral in (\ref{e.GK}) as a Riemann sum, comparison with (\ref{e.GK_disc}) shows that $\hat\sigma^2\,\kappa\to\sigma^2$ for $\kappa\to 0$.  We remark that for non-zero $\kappa$ the diffusion coefficient ${\hat{\sigma}}^2$ may differ from the diffusion coefficient $\sigma^2$ of the continuous system. 
Rescaling time to be measured in units of the discretization ``time step" $\kappa$, (\ref{e.SDE_map0}) can be rewritten as
\begin{align} 
\label{e.SDE_map}
dX=\Bigl(F(X)-\frac12 \kappa\, h(X)h'(X)\,\E[f_0^2] \Bigr)\, dt + \sqrt{\kappa}\, \hat\sigma h(X)\circ dW_t\, ,
\quad X(0)=\xi,
\end{align}
where $W_t$ is unit 1-dimensional Brownian motion on the rescaled time.

Comparing the limiting SDE of the discretized map (\ref{e.SDE_map}) and the limiting SDE for its associated continuous-time system (\ref{e.SDE_flow}), we see that they differ by an extra drift term in (\ref{e.SDE_map})
\begin{align}
E = -\frac12 \kappa \, h(X)h'(X)\,\E[f_0^2].
\label{e.extra}
\end{align}
Note that the additional drift term prohibits a Stratonovich interpretation of the noise and hence for finite $\Delta t$ the statistics of the map is different from the statistics of the original continuous-time system \eqref{e.ms.1}--\eqref{e.ms.2}, which we identify with the diffusive limit system (\ref{e.SDE_flow}). A discrepancy of this form was noted in \cite{GottwaldMelbourne13c} in the case of general maps where the time step $\Delta t$ (or rescaled time step $\kappa$) was implied. In the case where the fast dynamics of the discrete Euler map (\ref{e.map}) is {\em{i.i.d.}}, i.e.~$\hat \sigma^2=\E[f_0 ^2]$, the noise in the limiting SDE of the discretized map (\ref{e.SDE_map}) is of the It\^o type. This can be heuristically understood by realizing that if the time step $\kappa\gg \tau_{\rm{corr}}$ where $\tau_{\rm{corr}}$ is the decorrelation time of the fast continuous-time $y$-dynamics, the map is already as rough as the discrete approximation of the noise.\\

Although the additional drift term (\ref{e.extra}) is formally of order ${\mathcal{O}}(\Delta t)$, depending on the dynamical system under consideration, the extra term $E$ can be large and distort the statistical behavior leading to a marked difference between the numerically observed statistical behavior of the slow variable $X$ and the statistical behavior of the slow variable $x$ of the given continuous-time multi-scale system (\ref{e.ms.1})--(\ref{e.ms.2}) to be modelled. In Section~\ref{s.numerics} we provide such an example and show how an Euler discretization may produce erroneous statistical information.  In the next section we develop a relationship between the extra drift term (\ref{e.extra}) obtained in homogenization and backward error analysis. In particular, we will show that the extra drift term (\ref{e.extra}) generated by a first order numerical time-stepper is present in the backward error analysis and would be absent if the dynamics had been integrated with a higher order scheme instead (however, other terms are typically present in this case). This will show how the first order errors of an Euler-method directly translate into errors of the mean behavior.


\section{Backward error analysis}
\label{s.bea}

In this section we provide a backward error analysis to explain the presence of the extra term (\ref{e.extra}) in the homogenized discrete model (\ref{e.SDE_map}) compared to the homogenized continuous-time model (\ref{e.SDE_flow}).  We will see that the extra term arises from the use of the forward Euler scheme for constructing the discrete model (\ref{e.map}). 
Although the extra term (\ref{e.extra}) is of order $\mathcal{O}(\Delta t)$ and hence disappears  in the small step size limit, we stress that in the context of multi-scale problems one is often interested in step size regimes that are insensitive to fast dynamics and $\mathcal{O}(1)$ with respect to the slow dynamics. 

Backward error analysis \cite{HaLuWa06,LeimkuhlerReich} has been successfully employed to understand finite time step effects observed in numerical simulations.  
The truncation error of a numerical discretization of an ordinary differential equation (ODE) can be expanded as an asymptotic series in the step size $\Delta t$ with terms involving successively higher derivatives of the vector field. In backward error analysis, the terms of the truncated series are interpreted as a higher order approximation to another, perturbed vector field. 

\subsection{Lowest order modified equations}
We consider a generic differential equation
\begin{equation}\label{e.genODE}
	\dot{z} = \eff(z),
\end{equation}
the solution of which is to be approximated using a numerical method.  

To understand the qualitative behavior of the numerical solution for finite step size $\Delta t$,
one constructs a modified vector field as an asymptotic expansion
\begin{equation} \label{modODE}
	\dot{z} = \tilde{\eff}(z) = \eff(z) + \Delta t\,\eff_1(z) +  \Delta t^2\,\eff_2(z) + \cdots,
\end{equation}
where the terms $\eff_1$, $\eff_2$, etc.~are to be determined such that the solution to (\ref{modODE}) matches the expansion of a numerical method applied to (\ref{e.genODE}) to a higher order of accuracy.  The continuous-time solutions to the truncated modified differential equation (\ref{modODE}) approximate to higher order the numerical output than do those of the original differential equation (\ref{e.genODE}), allowing the modified equation to be used to interpret finite time step effects observed in the numerical time series.

The solution to (\ref{modODE}) is expanded in a Taylor series about $z(t)$ to give
\begin{equation}\label{modEqTaylor}
	z(t+\Delta t) = z(t) + \Delta t \, \tilde{\eff} + \frac{\Delta t^2}{2} \tilde{\eff}' \tilde{\eff} + \frac{\Delta t^3}{6} \left[ \tilde{\eff}''(\tilde{\eff},\tilde{\eff}) + \tilde{\eff}'\tilde{\eff}'\tilde{\eff} \right] + \mathcal{O} (\Delta t^4),
\end{equation}
where all terms on the right are evaluated at $z(t)$, and where $\tilde{\eff}'$ denotes the Jacobian matrix of $\tilde{\eff}$, $\tilde{\eff}''$ denotes its (symmetric) three-tensor of second partial derivatives, and $\tilde{\eff}''(\cdot,\cdot)$ denotes the contraction of this tensor with the two vector arguments.  Substituting the expansion (\ref{modODE}) into the above and gathering terms of like order yields
\begin{multline}\label{ModEqExp}
	z(t+\Delta t) = z(t) + \Delta t\, \eff + \Delta t^2 \left[ \eff_1 + \frac{1}{2} \eff'\eff \right]  + \\
	 \Delta t^3 \left[ \eff_2 + \frac{1}{2}(\eff'\eff_1 + \eff_1'\eff) 	+ \frac{1}{6} (\eff''(\eff,\eff) + \eff'\eff'\eff) \right] + \mathcal{O} (\Delta t^4).
\end{multline}
Next, one determines the functions $\eff_1$, $\eff_2$, etc.~to match the expansion of a numerical integrator to higher order.

Euler's method is given by $z_{n+1} = z_n + \Delta t\,\eff(z_n)$.  This formula is consistent with (\ref{ModEqExp}) up to terms of $\mathcal{O}(\Delta t^2)$, and is consequently a first order approximation to (\ref{e.genODE}).  However by choosing
\[
	\eff_1 = - \frac{1}{2} \eff'\eff
\]
in (\ref{modODE}),
one finds that Euler's method agrees with (\ref{ModEqExp}) up to terms of $\mathcal{O}(\Delta t^3)$.  Consequently, while Euler's method is a first order approximation of (\ref{e.genODE}), it is a \emph{second order approximation} to the modified differential equation 
\begin{align}
	\dot{z} = \eff - \frac{\Delta t}{2} \eff'\eff.
\label{e.Eulermod}
\end{align}
This process may be repeated to derive higher order corrections ($\eff_2$, $\eff_3$, etc.) in the modified equation. The asymptotic expansion generally does not converge for fixed $\Delta t$, but may be optimally truncated \cite{HaLuWa06}. Although for general systems there is no guarantee that the lowest order corrections will have the most significant impact on the observed statistics, in our numerical experiments this does appear to be the case.

Next consider the second order Runge-Kutta method (i.e.~Heun's method)
\begin{align}
	z_{n+1} = z_n + \frac{\Delta t}{2} \left[ \eff(z_n) + \eff( z_n + \Delta t \eff(z_n) ) \right].
\label{e.RKex}
\end{align}
Expanding the right-hand side about $z_n$ gives
\[
	z_{n+1} = z_n + \frac{\Delta t}{2} \left[ 2\eff(z_n)+ \Delta t \eff'(z_n) \eff(z_n) + \frac{\Delta t^2}{2} \eff''(\eff,\eff)(z_n)  + \mathcal{O}(\Delta t^3) \right].
\]
By choosing
\[
	\eff_1 = 0, \qquad 	\eff_2 = \frac{1}{12} \eff''(\eff,\eff) -\frac{1}{6} \eff'\eff'\eff 
\]
in (\ref{modODE}) we make this formula agree with (\ref{ModEqExp}) up to terms of $\mathcal{O}(\Delta t^4)$ so the modified equation associated with the second order Runge-Kutta method (\ref{e.RKex}) is (up to terms of $\mathcal{O}(\Delta t^2)$)
\begin{align}
	\dot{z} = \eff(z) + \frac{\Delta t^2}{12} ( \eff''(\eff,\eff) - 2 \eff'\eff'\eff ),
\label{e.RKmod}	
\end{align}
and the second order Runge-Kutta method (\ref{e.RKex}) applied to the ODE (\ref{e.genODE}) is a third order approximation to (\ref{e.RKmod}).\\

Note that the modified equation (\ref{e.Eulermod}) suggests a correction to Euler's method to eliminate the second order term in (\ref{ModEqExp}).  Specifically, one can apply Euler's method to the differential equation with corrected vector field
\[
	\dot{z} = \eff(z) + \frac{\Delta t}{2} \eff'(z) \eff(z).
\]
Doing so yields the second order Taylor method
\begin{equation}\label{e.Taylor}
	z_{n+1} = z_n + \Delta t \, \eff(z_n) + \frac{\Delta t^2}{2} \eff'(z_n) \eff(z_n),
\end{equation}
which can be efficiently implemented using a finite difference approximation in the last term:
\[
	\eff'(z_n)\eff(z_n) \approx \frac{1}{\tau} \left( \eff(z_n + \tau \eff(z_n)) - \eff(z_n) \right), \quad \tau = \sqrt{\epsilon_m},
\]
with machine precision $\epsilon_m$. Matching (\ref{e.Taylor}) with (\ref{ModEqExp}) shows the Taylor method has modified equation expansion (up to terms of $\mathcal{O}(\Delta t^2)$)
\begin{equation}\label{e.modEqTaylor}
	\dot{z} = \eff(z) - \frac{\Delta t^2}{6} ( \eff''(\eff,\eff) + \eff'\eff'\eff),
\end{equation}
and is hence is a second order scheme for the original ODE (\ref{e.genODE}). This will turn out to be advantageous for problems of the form considered here.

\subsection{Homogenization of modified equations}
In the limit $\varepsilon\to 0$, the time step scales as $\Delta t = \kapfactor\,\varepsilon^2$, where $\kapfactor>0$ is fixed and small.  As we will take this limit to homogenize the modified equation, we use $\kapfactor = \Delta t / \varepsilon^2$ as our expansion parameter for the backward error analysis.  

The modified vector fields $\eff_1$, $\eff_2$, etc.~each in turn can be expressed as expansions in $\varepsilon$. Upon homogenization, the lowest order term of $\mathcal{O} (1/\varepsilon)$ contributes to the diffusion and the $\mathcal{O}(1)$ term contributes to the drift. Terms of higher order in $\varepsilon$ vanish in the homogenization limit $\varepsilon\to 0$. 

Substituting the vector fields of the deterministic multi-scale system (\ref{e.ms.1})--(\ref{e.ms.2}) into the modified equation for a first order Euler discretization \eqref{e.Eulermod} with $z=(x,y)$ we obtain 
\begin{align}
\dotx = &\frac{1}{\varepsilon}\, h(\x)f_0(\y)+ f(\x,\y)  \label{e.BE_x}\\
&- \frac{\kapfactor}{2}\left(   \frac{1}{\varepsilon}\, h(x) f_0'(y)\, \tilde{g}(y) + 
h(x)  h'(x)  f_0^2(y) + \partial_y f(x,y) \tilde{g}(y) + \mathcal{O}(\varepsilon)\right ) 
+\mathcal{O}(\kapfactor^2) \nonumber\\
	\dot{y} = &\frac{1}{\varepsilon^2}\tilde{g}(y), \qquad \tilde{g} =  g(y) + \frac{\kapfactor}{2K}  g'(y)g(y) + \mathcal{O}((\kapfactor/K)^2).
\end{align}
In the remainder of the Section we substitute the vector field $g(y)$ for the vector field $\tilde{g}(y)$ of the fast modified equation, and similarly substitute $\varphi^t$ for $\tilde\varphi^t$, the flow of $\tilde{g}$, which is again independent of $\varepsilon$ since $\kapfactor$ and $K$ are fixed. This is admissible if the numerical scheme for the fast dynamics is sufficiently accurate to resolve its statistical behaviour and in particular, the autocorrelation function; given the relation $\delta t = \Delta t/K$ and the scaling $\Delta t =\kapfactor\, \varepsilon^2$, it is clear that $\kapfactor/K$ should be chosen small enough to accurately approximate autocorrelation functions of the unscaled chaotic system with vector field $g(y)$, to allow for the substitution of $\varphi^t$ for $\tilde\varphi^t$.  

We now show that the associated stochastic limit system of the continuous-time modified equation of the Euler method (\ref{e.Eulermod}) is the same as that of the discrete Euler discretization 
(\ref{e.map}). 
The homogenized dynamics, approximating the dynamics of (\ref{e.BE_x}) on time scales of order ${\mathcal{O}}(1)$, is given up to $\mathcal{O}(\kapfactor)$ by
\begin{align}
dX=F(X)\,dt + \sigma h(X)\circ dW_t, \quad X(0)=\xi.
\label{e.BE_homo}
\end{align}
where $F(X)$ is the expectation with respect to $y\sim\mu$ of the terms of $O(1)$ in $\varepsilon$ in \eqref{e.BE_x}.  Noting that ${g}(y) = \eps^{2}\doty$, we observe that the terms $f_0^\prime(y) {g}(y)$ and $\partial_y f(x,y) {g}(y)$ in the slow modified equation (\ref{e.BE_x}) can be written as a total derivative (taking $x$ constant up to terms of $\O(\varepsilon)$ on the homogenization time scale).  Consequently, these terms vanish in expectation as in, for example,  
\begin{equation}\label{e.totalD}
	\E\left[ \partial_y f(X,y){g}(y) \right] = \E\left[ \frac{d}{dt} f(X,y(t))\big|_{X \mathrm{fixed}} \right] = 0,
\end{equation}
and we are left with the drift term
\begin{align}
F(X)&=\E[f(X,y)] - \frac{\kapfactor}{2} h(X) h'(X) \,\E[f_0^2] 
\nonumber
\\
&= \E[f(X,y)] + E.
\label{e.BE_homoF}
\end{align}
The diffusion coefficient is given by
\begin{multline}
\frac{1}{2}\sigma^2 = 
\int_0^\infty 
\E\left[ f_0(y)f_0(\tilde\varphi^ty) \right] dt
\\
+ \frac{\kapfactor}{2}h(X) \int_0^\infty \left(
\E\left[ f_0(y)\,(\frac{d}{dt}f_0)(\tilde\varphi^ty)  \right]
+
\E\left[ f_0(\tilde\varphi^ty)\,(\frac{d}{dt}f_0)(y) \right] 
\right)\, dt\, .
\label{e.GK_BE}
\end{multline}
Using ergodicity of the fast dynamics the spatial average in the second integral can be expressed as a time-average; partial integration then can be used to show that the second integral sums to zero. The Green-Kubo formula (\ref{e.GK_disc}) in the limit of $\varepsilon \to 0$ is recovered provided we substitute $\varphi^t$ for $\tilde\varphi^t$ in (\ref{e.GK_BE}), which we argued above is admissible provided $\kappa/K\ll 1$. 
Hence in the limit $\varepsilon\to 0$ we recover the homogenized equation (\ref{e.SDE_map}) for the forward Euler map (\ref{e.map}), and the long term statistics of the Euler discretization captures well the statistics of its associated continuous-time modified equation.


For a second order method, such as the Runge-Kutta method (\ref{e.RKex}) or second order Taylor method (\ref{e.Taylor}), the additional drift term $E$ is absent from the modified equation.  Nevertheless there are terms of $\mathcal{O}(\kapfactor^2)$ that could potentially influence the homogenized limit.  
When the Runge-Kutta method (\ref{e.RKex}) is applied to the deterministic multi-scale system (\ref{e.ms.1})--(\ref{e.ms.2}), the lowest order modified equation (\ref{e.RKmod}) yields, for the slow variable: 
\begin{multline}
\dotx = \frac{1}{\varepsilon}\, h(\x)f_0(\y)+ f(\x,\y)  \\
+ \frac{\kapfactor^2}{12} \left( \frac{1}{\varepsilon}\left[h(x)\left(f_0''(y)g^2(y) - 2f_0'(y)g'(y)g(y)\right) \right] \right. \\
\left. \phantom{\frac{\kapfactor^2}{12}}+ \left[ \partial_{yy}f(x,y) g^2(y) - 2 \partial_yf(x,y) g'(y)g(y)\right] + \mathcal{O}(\varepsilon) \right)
+\mathcal{O}(\kapfactor^3).
\label{e.BE_x_rk}
\end{multline}
The homogenized drift term becomes
\[
	F(X) = \mathbb{E} \left[ f(X,y) + \frac{\kapfactor^2}{12} \left( \partial_{yy} f(X,y)g^2(y) - 2 \partial_y f(X,y) g'(y)g(y) \right) \right] +\mathcal{O}(\kapfactor^3),
\]
which implies an $\mathcal{O}(\kapfactor^2)$ bias in the drift of the slow variables in the limit $\varepsilon\to 0$. The homogenized diffusion parameter becomes
\[
	\frac{1}{2}\sigma^2 = h(X) \int_0^\infty \C \left[ f_0(y) + \frac{\kapfactor^2}{12} \left( f_0''(y)g^2(y) - 2 f_0'(y)g'(y)g(y) \right) \right](t) \, dt +\mathcal{O}(\kapfactor^3).
\]
Here, too, an $\mathcal{O}(\kapfactor^2)$ bias occurs as $\varepsilon\to 0$.

Finally, when the Taylor method (\ref{e.Taylor}) is applied to the deterministic multi-scale system (\ref{e.ms.1})--(\ref{e.ms.2}), the lowest order modified equation (\ref{e.modEqTaylor}) yields, for the slow variable: 
\begin{equation}
\dotx = \frac{1}{\varepsilon}\, h(\x)f_0(\y)+ f(\x,\y)  
+ \frac{\kapfactor^2}{6}\left( 
\frac{1}{\varepsilon} v_1^{\mathrm{diff}}(x,y) + v_1^{\mathrm{drift}}(x,y) + \mathcal{O}(\varepsilon) \right) + \mathcal{O}(\kapfactor^3), \label{e.BE_x_taylor}
\end{equation}
where
\[
 v_1^{\mathrm{diff}}(x,y) = h(x)f_0''(y)g^2(y) + h(x)f_0'(y)g'(y)g(y),
\]
and 
\[
 v_1^{\mathrm{drift}}(x,y)=  3 h(x)'h(x)f_0^\prime (y)f_0(y)g(y) + \partial_{yy}f(x,y)g^2(y) + \partial_y f(x,y)g'(y)g(y).
\]
Noting that $g(y) = \eps^{2}\doty$, we observe that all terms in the drift perturbation can be written as total derivatives with $x$ fixed and vanish in expectation (cf.~\eqref{e.totalD}):
\[
	\E[v_1^{\mathrm{drift}}(X,y)] = 
	\E\left[ 3 h'(X)h(X)\frac{d}{dt}\left(\frac{f_0^2(y)}{2}\right) 
	+ \frac{d}{dt} \left( \partial_y  f(X,y)g(y) \right)_{X \mathrm{fixed}} \right] = 0.
\]
Consequently, the Taylor method has \emph{no bias in drift to $\mathcal{O}(\Delta t^3)$}, and we expect the drift to be simply given by
\begin{equation}\label{e.drift_Taylor}
	F(X) = \mathbb{E} \left[ f(X,y) \right].
\end{equation}
The diffusion term is also a total derivative: 
\[
	 v_1^{\mathrm{diff}}(x,y) = \frac{1}{\varepsilon} h(X) \frac{d}{dt} \left( f_0'(y)g(y) \right),
\]
and the diffusion parameter $\sigma$ is of the form 
\[
	\frac{1}{2}\sigma^2 = h(X) \int_0^\infty \left( \C[f_0](t) - \frac{\kapfactor^2}{3} \C\left[\frac{d}{dt} f_0\right](t)\right) \, dt + \mathcal{O}(\kapfactor^3).
\]
The term does not vanish; this is a correlation function, not just an expectation.

In summary, we note that the additional drift term $E$ of order $\mathcal{O}(\kapfactor)$ in (\ref{e.BE_homoF}) is absent in the modified equation of a numerical method which is at least second order, such as the Runge-Kutta method (\ref{e.RKex}) or the second order Taylor method (\ref{e.Taylor}). For a second order time-stepping method the homogenized modified equation therefore agrees with the homogenized equation of the full multi-scale system (\ref{e.SDE_flow}) up to $\mathcal{O}(\kapfactor^2)$. However, a second order scheme will generally also have additional corrections to the drift and diffusion which might be of the same magnitude for finite $\kapfactor$ as those corresponding to the continuous-time multi-scale system under consideration.  Remarkably, the second order Taylor method (\ref{e.Taylor}) does not have bias of $\mathcal{O}(\kapfactor^2)$ in the drift. This can be traced to the fact that the scheme exactly agrees with the second order Taylor expansion of the error, and consequently its higher order terms are exact differentials. 

\section{Numerical demonstration}
\label{s.numerics}
We now demonstrate that the additional terms in the backward error analysis may lead to significant bias in the probability density estimation for both the first order forward Euler scheme and, to a lesser degree, the second order Heun's method. In particular we show that the numerical methods converge for $\eps\to 0$ with $\Delta t = \kapfactor \, \eps^2$ to the homogenized limits of their respective modified equations, which are different from the long-time statistical limit of the original deterministic multi-scale system.\\ We consider the deterministic multi-scale system (\ref{e.ms.1})--(\ref{e.ms.2}) with
\[
f_0(y)=a y, \quad h(x)=\sqrt{x}, \quad
f(x,y)=b(c-x)y^2,
\]
so the slow dynamics is described by the continuous-time system
\begin{align}
\dotx &= \frac{1}{\varepsilon} a \sqrt{\x} y + b (c-\x)y^2\, .
\label{e.CIR_ms}
\end{align}
We choose here $a = 0.1$, $b = 0.005$ and $c = 0.75$. The slow dynamics is driven by $y=\zeta_2 + \zeta_3$ generated by a fast chaotic R\"ossler system
\begin{align}
\eps^2\dot \zeta_1 &= -\zeta_2-\zeta_3, \label{e.CIRa}\\
\eps^2\dot \zeta_2 &= \zeta_1+r\zeta_2, \label{e.CIRb}\\
\eps^2\dot \zeta_3 &= s +(\zeta_1-u)\zeta_3\, , \label{e.CIRc}
\end{align}
with $r=s=0.25$ and $u=7$. We will compare the results of a numerical integration of this deterministic multi-scale system using first and second order discretization methods to results from the associated limiting homogenized SDE of this system, describing the long-time statistical behavior.\\
The diffusive limiting equation of the multi-scale dynamical system (\ref{e.CIR_ms})--(\ref{e.CIRc}) 
can be obtained via the homogenization techniques presented in Section \ref{s.ms} and is given by the Cox-Ingersoll-Ross (CIR) model~\cite{CIR85,CIR85b} 
\begin{align} 
dX = \sigma a \sqrt{X} \, dW_t +2\alpha b (\beta-X)\, dt,
\label{e.CIR}
\end{align}
where $W_t$ is unit $1$-dimensional Brownian motion. The parameters are:
\begin{align}
\alpha&=\frac12\,\E[ y^2], \label{e.para} \\
\sigma^2&=2\int_0^\infty \E[(\varphi^ty)y]\, dt, \label{e.sig_cont}\\
\beta&=c+\frac{\sigma^2 a^2}{8\alpha b}. \label{e.bet_cont}
\end{align}
To approximate $\alpha$, an ensemble simulation of the (unscaled) R\"ossler system was carried out using a 1000-member ensemble on a time interval $t\in[0,3.2\times 10^4]$ with initial conditions drawn approximately from $\mu$ (see below). 
We obtain $\alpha = 28.4\pm 0.1$. To estimate $\sigma^2$ we solve $w_{n+1}=w_n+\varepsilon\, \Delta t\, y_n$ where $\Delta t=\kapfactor\,\eps^2$ and $\kapfactor=0.5$. Then for a time series of length $N=\lfloor 1/\varepsilon^2 \rfloor$, $w_N \sim \tfrac{\Delta t}{\sqrt{N}}\sum_{j=1}^{N-1} y_j$ is approximate Brownian motion with variance $ \V[w_N] = \sigma^2\, N \Delta t$. In this way the diffusivity is estimated as $\sigma^2 \approx 0.140 \pm 0.002$.

The Cox-Ingersoll-Ross (CIR) model (\ref{e.CIR}) has the closed form solution 
\begin{align} 
\label{e.CIR_sol}
\frac{X(t)}{c(t)} \sim H(t),\quad c(t)=\frac{\sigma^2a^2}{8\alpha b}(1-e^{-2\alpha b t}),
\end{align}
where $H(t)$ is a noncentral $\chi$-squared distribution with $8\alpha\beta b/(a^2 \sigma^2)$
degrees of freedom and noncentrality parameter 
$c(t)^{-1}e^{-2\alpha bt}X(0)$. \\

To numerically integrate the two-scale system (\ref{e.CIR_ms})--(\ref{e.CIRc}) we use a multiple time-stepping approach \cite{LeimkuhlerReich}, with a step size $\Delta t=\kapfactor\, \eps^2$ for (\ref{e.CIR_ms}) and step size $\delta t = \Delta t/K$ for the fast subsystem (\ref{e.CIRa})--(\ref{e.CIRc}). For our illustration we choose successively $\varepsilon \in \{ 0.05, 0.025, 0.0125, 0.00625\}$ and integrate over the interval $t\in[0,2.5]$ using $\kapfactor = 0.5$ and $K=50$.  For this scaling of time step the fast dynamics (\ref{e.CIRa})--(\ref{e.CIRc}) is well resolved but is not solved with increasing accuracy in the limit $\eps\to 0$. The probability density function (pdf) of $x(t)$ is estimated using an ensemble with 160000 members.  Each member starts from $x(0)=1$ but observes a distinct time series $y(t)$.
Each member initial condition $y(0)$ is drawn from the invariant measure $\mu$ by letting a randomly drawn initial condition relax onto the attractor over a transient time of length $25$. 
We initially compare two different numerical discretizations of the multi-scale system (\ref{e.CIR_ms})--(\ref{e.CIRc}). Applied to the generic differential equation
\[
	\dot{x}(t) = \eff (x(t),y(t)),
\]
these methods are: the forward Euler method
\begin{equation} \label{e.fe}
	x_{n+1} = x_n + \Delta t \, \eff (x_n, y_n), 
\end{equation}
and the second order Runge-Kutta method (Heun's method)
\begin{equation}\label{e.etr}
	x_{n+1} = x_n + \frac{\Delta t}{2} \left[ \eff(x_n,y_n) + \eff\left(x_n + \Delta t \eff(x_n,y_n), y_{n+1}\right) \right].
\end{equation}
In the above equations $y_n$ denotes the approximation to $y(t_n)$ obtained from $nK$ steps of size $\delta t$.\\ When the slow dynamics (\ref{e.CIR_ms}) is discretized using a forward Euler method (\ref{e.fe}) with time step $\Delta t$ we obtain the map 
\begin{align}
x_{n+1} &= x_n + \Delta t \frac{1}{\varepsilon} \, a\sqrt{x_n} y_n +\Delta t \, b (c-x_n) y_n^2.
\label{e.CIR_disc}
\end{align}

We compare the following probability density functions at time $t=2.5$:
\begin{itemize}
\item{} [MS1] The empirical pdf of the multi-scale system (\ref{e.CIR_ms})--(\ref{e.CIRc}) computed using the forward Euler scheme (\ref{e.fe}).
\item{} [MS2] The empirical pdf of the multi-scale system (\ref{e.CIR_ms})--(\ref{e.CIRc}) computed using the second order Runge-Kutta scheme (\ref{e.etr}).
\item{} [HMC] The exact pdf (\ref{e.CIR_sol}) of the limiting homogenized stochastic CIR model (\ref{e.CIR}) with parameters $\alpha$ given by (\ref{e.para}) and $\sigma$ and $\beta$ given by (\ref{e.sig_cont})--(\ref{e.bet_cont}) associated with the continuous-time model  (\ref{e.CIR_ms})--(\ref{e.CIRc}).
\end{itemize}
Empirical pdfs of the slow variable $x$ are computed using Matlab's histogram counter with bin size $\Delta x = 0.005$.

\begin{figure}[!h]
\begin{center}
\includegraphics[width=0.49\textwidth]{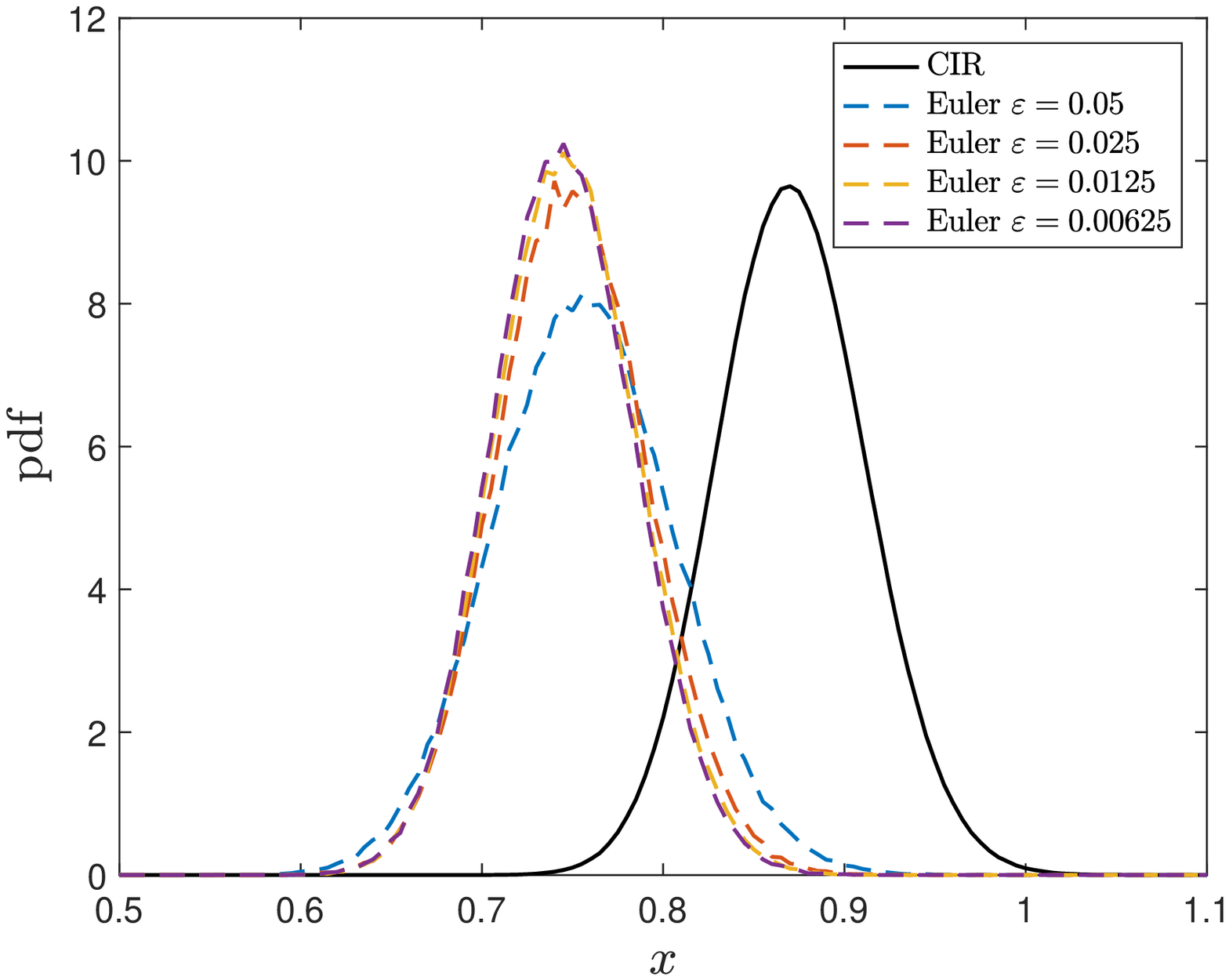}
\includegraphics[width=0.49\textwidth]{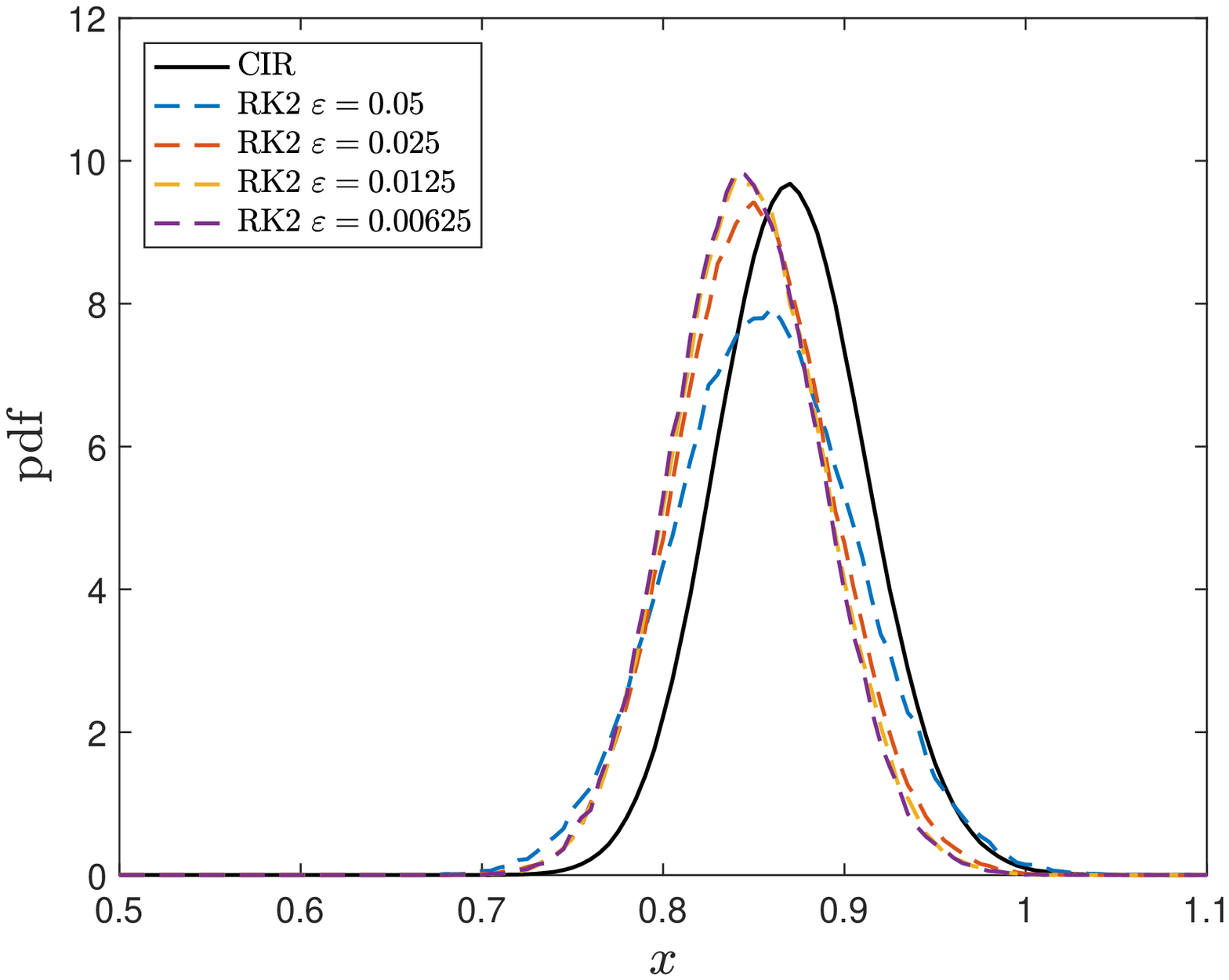}
\caption{Comparison of the pdfs of the numerically approximated multi-scale system  (\ref{e.CIR_ms})--(\ref{e.CIRc}) with the exact density of the CIR model (\ref{e.CIR_sol}) at time $t=2.5$ for the forward Euler scheme (\ref{e.fe}) (left) and the second order Runge-Kutta scheme (\ref{e.etr}) (right). In each plot the solid black line indicates the exact pdf associated with the homogenized CIR model (\ref{e.CIR}) 
for the actual time-continuous multi-scale system (\ref{e.CIR_ms})--(\ref{e.CIRc}) (i.e. with parameters (\ref{e.para})--(\ref{e.bet_cont})). The dashed lines indicate the empirical probability density function for the numerical simulations of the continuous-time multi-scale system (\ref{e.CIR_ms})--\eqref{e.CIRc}
for $\varepsilon=0.05$ (blue), $0.025$ (red), $0.0125$ (yellow), and $0.00625$ (purple). 
\label{fig:compCIR}}
\end{center}
\end{figure}

The rigorous homogenization results presented in Section~\ref{s.ms} assert that the long-time statistics of the full deterministic multi-scale system is described by the pdf [HMC]. 
Figure \ref{fig:compCIR} shows convergence of the empirical pdfs of the numerically approximated multi-scale system, obtained with the forward Euler ([MS1], left frame) and Runge-Kutta ([MS2], right frame) methods for decreasing $\varepsilon$ (dashed lines, with $\varepsilon=0.05$ in blue, $\varepsilon=0.025$ in red,  $\varepsilon=0.0125$ in yellow, and $\varepsilon=0.00625$ in purple). The exact pdf [HMC] is also indicated in each frame (solid black line).  Both numerically computed pdfs appear to converge, in the limit $\varepsilon \to 0$, $\Delta t =\kapfactor\,\varepsilon^2$, to a density with the wrong mean.   The pdf [MS2] of the RK method is significantly closer to [HMC] than is [MS1], but bias is nevertheless clearly present.

\begin{figure}[!h]
\begin{center}
\includegraphics[width=0.49\textwidth]{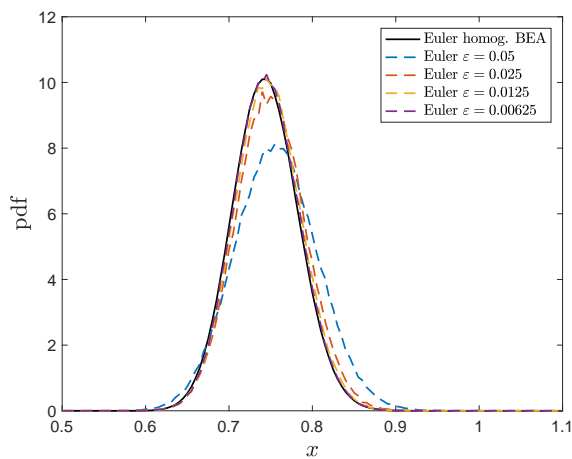}
\caption{Comparison of the pdfs of the numerically approximated multi-scale system (\ref{e.CIR_ms})--(\ref{e.CIRc}) with the exact density of the homogenized modified equation model at time $t=2.5$ for the forward Euler scheme (\ref{e.fe}). The solid black line indicates the exact pdf associated with the homogenized modified equation (\ref{e.CIR})
of the Euler scheme (i.e. with parameters (\ref{e.para_disc})-(\ref{e.bet_disc})). 
The dashed lines indicate the empirical probability density function computed with the forward Euler approximation (\ref{e.CIR_disc}) 
applied to the continuous-time multi-scale system (\ref{e.CIR_ms})--(\ref{e.CIRc}) for $\varepsilon=0.05$ (blue), $0.025$ (red), $0.0125$ (yellow), and $0.00625$ (purple). 
\label{fig:compBEA}}
\end{center}
\end{figure}

For the first order forward Euler discretization (\ref{e.CIR_disc})
the homogenized SDE (\ref{e.BE_homo}) describing the long-time behavior of the slow motion is also given by the CIR model (\ref{e.CIR}), but now with parameters
\begin{align}
\alpha&=\frac12\,\E[ y^2], \label{e.para_disc} \\
\hat \sigma^2&= \E[ y^2]+2\sum_{n=1}^\infty \E[ (\Phi^n\,y)y]=\lim_{n\to\infty}n^{-1} \E[ (\sum_{j=0}^{n-1}\Phi^j\,y)^2], \label{e.sig_disc}\\
\beta&=c+\frac{\hat\sigma^2\Delta t \,a^2}{8\alpha b} - \frac{\Delta t a^2 }{4b}, \label{e.bet_disc}
\end{align}
and $\hat\sigma^2 \Delta t \approx \sigma^2$. 
Note that the only difference in the homogenized CIR systems associated with the continuous-time multi-scale system (\ref{e.ms.1})--(\ref{e.ms.2}) and the discrete Euler map (\ref{e.CIR_disc}) is in the parameter $\beta$ (cf. (\ref{e.bet_cont}) and (\ref{e.bet_disc})). For $a^2/b\gg 1$ and $\sigma^2/4\alpha \ll 1$ this difference is large and may cause significant discrepancy between the statistics of the continuous-time multi-scale system \eqref{e.CIR_ms}--\eqref{e.CIRc} and its first order Euler discretization \eqref{e.CIR_disc}. The latter condition, $\sigma^2/4\alpha=\int_0^\infty \E[(\varphi^ty)y]\, dt/\E[  y^2]  \ll 1$ puts a requirement on the decay of the fast dynamics and states that the fast dynamics should be far from {\em{i.i.d.}}~with $\sigma^2=2 \, \E[ y^2 ]$. This requirement is satisfied for the R\"ossler system \eqref{e.CIRa}--\eqref{e.CIRc} with the parameters $r=s=0.25$ and $u=7$. For these parameters the autocorrelation function has a slow decay and $\E[y^2]\gg \sigma^2/2$. 
For comparison, we introduce an additional pdf:
\begin{itemize}
\item{}  [HMD] The exact pdf (\ref{e.CIR_sol}) of the limiting homogenized stochastic CIR model (\ref{e.CIR}) with parameters $\alpha$ given by (\ref{e.para}) and $\sigma$ and $\beta$ given by (\ref{e.sig_disc})--(\ref{e.bet_disc}) associated with the discrete Euler model  (\ref{e.CIR_disc}).
\end{itemize}
The backward error analysis presented in Section~\ref{s.bea} predicts that the empirical pdf [MS1] of the Euler method may be better approximated by the pdf [HMD], derived by homogenizing the 
modified equation of the Euler method.  
Figure~\ref{fig:compBEA} confirms this prediction, showing that as $\varepsilon\to 0$, the statistical behavior of the discrete Euler scheme is well described by the pdf of its associated homogenized stochastic CIR model.
The extra drift term $E$ in the homogenized discrete model leads to an error of $16\%$ in the mean of the pdf [HMD] with respect to the mean of the pdf of the original continuous time multi-scale system  (\ref{e.CIR_ms})--(\ref{e.CIRc}) to be modelled [HMC].

We remark that for fast chaotic dynamics with rapidly decaying autocorrelation function such as the Lorenz $63$ system with the classical parameters, we have $\E[y^2]\approx \sigma^2/2$.  The homogenized equation of the full multi-scale dynamics and its first order Euler discretization will be close (cf.~(\ref{e.bet_cont}) and (\ref{e.bet_disc})), and a first order discretization would be sufficient to capture the long-time statistics of the slow dynamics.

\begin{figure}[!h]
\begin{center}
\includegraphics[width=0.49\textwidth]{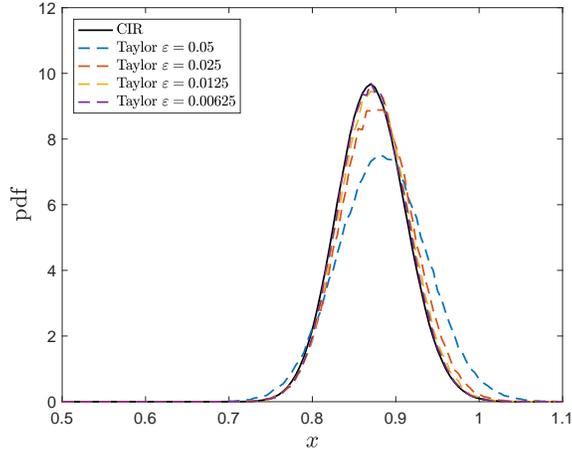}
\caption{Comparison of the pdfs of the numerically approximated multi-scale system  (\ref{e.CIR_ms})--(\ref{e.CIRc}) with the exact density of the CIR model (\ref{e.CIR_sol}) at time $t=2.5$ for the second order Taylor method (\ref{e.Taylor}). The solid black line indicates the exact pdf associated with the homogenized CIR model (\ref{e.CIR}) 
for the actual time-continuous multi-scale system (\ref{e.CIR_ms})--(\ref{e.CIRc})  (i.e. with parameters (\ref{e.para})--(\ref{e.bet_cont})). The dashed lines indicate the empirical probability density function computed with the second order Taylor method (\ref{e.Taylor}) applied to the continuous-time multi-scale system (\ref{e.CIR_ms})--(\ref{e.CIRc})
for $\varepsilon=0.05$ (blue), $0.025$ (red), $0.0125$ (yellow), and $0.00625$ (purple). 
\label{fig:compCorrEuler}}
\end{center}
\end{figure}

Finally, the backward error analysis of the second order Taylor method (\ref{e.Taylor}) indicates there is no error in the drift (\ref{e.drift_Taylor}) to $\mathcal{O}(\kapfactor^3)$ for this method.  Indeed, Figure \ref{fig:compCorrEuler} confirms this result, illustrating that the empirical pdf of the Taylor method closely matches that of the pdf [HMC] as $\varepsilon \to 0$.  In fact both the mean and the variance of the distribution closely match that of the pdf [HMC], suggesting that errors in the diffusion parameter are also small for this parameter regime.  The excellent approximation of the drift makes the Taylor method an attractive alternative for multi-scale problems with stochastic limit behavior.

\section{Summary}
\label{s.summary}

To summarize, using backward error analysis we have demonstrated that the extraneous drift  term (\ref{e.extra}) that arises in homogenization of the discrete map (\ref{e.map}) compared to homogenization of the flow (\ref{e.ms.1})--(\ref{e.ms.2}) can be traced to the relation between the map (\ref{e.map}) and a forward Euler discretization of (\ref{e.ms.1})--(\ref{e.ms.2}) with large step size $\Delta t=1$. In particular, we relate this drift term which appears in the literature of homogenization for discrete time systems \cite{GottwaldMelbourne13c} and which neither corresponds to an It\^o nor to a Stratonovich interpretation of the SDE to discretization errors of first order schemes using backward error analysis. We have shown that the local first order errors contribute to a well-defined drift error, leading to potentially strong bias, in the long-time statistical behavior. The accumulated local error, as quantified by the backward error analysis, was shown to account for the long-time statistical error of the discretization scheme as provided by homogenization theory.  

We further quantified the requirement for a dynamical system such that its Euler scheme discretization reliably recovers its long-time statistical behavior. In particular, we found that for sufficiently rapidly decaying fast dynamics an Euler scheme is sufficient. On the contrary, the failure of first order discretization methods to capture the statistics of the full continuous-time multi-scale system was shown to be exacerbated if the fast dynamics exhibits slow decay of correlations. We remark that the slow decay of correlation is not hampering the validity of the homogenized limit system and the validity of the underlying functional central limit theorem which is assured solely by requiring $\varepsilon \ll 1$. The difference is entirely given by the failure to match the limiting homogenized SDE of the first order discretization with the limiting homogenized SDE associated with the original time-continuous multi-scale system.

Here we discussed deterministic skew product systems of the form (\ref{e.ms.1})--(\ref{e.ms.2}). In order to obtain a stochastic homogenized equation for the slow dynamics, the fast dynamics is required to support an ergodic invariant measure and generate an integrable autocorrelation function of $f_0$ (cf.~the Green-Kubo formula (\ref{e.GK})). Hence the conclusions drawn here for the deterministic setting remain valid in the case when the fast dynamics is stochastic. For a large class of stochastic ordinary differential equations stochastic integrators were constructed which accurately approximate the invariant measure \cite{AbdulleEtAl14}. Their construction also uses the framework of modified equations and the analysis is not restricted to the multi-scale setting. It would be interesting to compare the higher order methods developed there in the multi-scale setting considered here. This is a topic for further research.\\ 

In this article, we have examined the limit $\varepsilon\to 0$, $\Delta t \sim \O({\varepsilon^2})$ in multi-scale systems (\ref{e.ms.1})--(\ref{e.ms.2}) which approach a rigorous SDE limit under homogenization.  In this limit, with constant $\kappa=\Delta t/\varepsilon^2>0$, the fast motion of the system remains resolved, but \emph{does not converge} as $\Delta t\to 0$.  The limit is achieved by choosing $K$ to be sufficiently large and $\kappa/K$ to be sufficiently small. Our numerical experiments demonstrate significant bias in the pdf of the slow variables in the limit $\Delta t \to 0$ under this scaling.  This bias may be mitigated by reducing the ratio $\kappa=\Delta t /\varepsilon^2$ in the numerical experiments as the bias is multiplied by $\kappa$ (cf. \ref{e.extra}).  However, full convergence requires $\kappa=\Delta t/\varepsilon^2\to 0$ as $\varepsilon\to 0$.\\ We assumed throughout that the fast dynamics is numerically sufficiently resolved such that the statistical properties, e.g.~the auto-correlation structure, is sufficiently reproduced. If this were not the case, errors arising from the flow map associated with the fast modified equation enter the Green-Kubo formulae (see, for example, (\ref{e.GK_BE}) for the Euler method) implying errors in the diffusion coefficient in addition to the bias error.

We remark that special numerical methods are often specifically tailored to multi-scale problems to accommodate time steps that are large with respect to the fast time scale. It is precisely in this regime that statistical bias may occur. The implication of this for numerical integration of multi-scale systems is that, to avoid statistical bias, it may be important to use a higher order method for the slow variables. The second order Taylor method (\ref{e.Taylor}) offers an interesting alternative here, as it may be efficiently implemented and is unbiased with respect to the drift up to $\mathcal{O}(\Delta t^3)$.
To avoid statistical bias altogether, one might want to solve the actual limiting SDE instead of the deterministic multi-scale system provided that $\epsilon$ is sufficiently small to allow for the central limit theorem to hold.

\section*{Acknowledgments}
We gratefully acknowledge funding through the University of Sydney -- University of Utrecht Partnership Collaboration Award.

\bibliographystyle{siamplain}

\begin{thebibliography}{10}

\bibitem{AbdulleEtAl14}
{\sc A.~Abdulle, G.~Vilmart, and K.~C. Zygalakis}, {\em High order numerical
  approximation of the invariant measure of ergodic sdes}, SIAM Journal on
  Numerical Analysis, 52 (2014), pp.~1600--1622.

\bibitem{AbKoMa03}
{\sc R.~V. Abramov, G.~Kova{\v{c}}i{\v{c}}, and A.~J. Majda}, {\em Hamiltonian
  structure and statistically relevant conserved quantities for the truncated
  burgers-hopf equation}, Communications on pure and applied mathematics, 56
  (2003), pp.~1--46.

\bibitem{ArielEtAl13}
{\sc G.~Ariel, B.~Engquist, S.~Kim, Y.~Lee, and R.~Tsai}, {\em A multiscale
  method for highly oscillatory dynamical systems using a {P}oincar\'e map type
  technique}, J. Sci. Comput., 54 (2013), pp.~247--268.

\bibitem{CIR85}
{\sc J.~C. Cox, J.~E. Ingersoll, Jr., and S.~A. Ross}, {\em An intertemporal
  general equilibrium model of asset prices}, Econometrica, 53 (1985),
  pp.~363--384.

\bibitem{CIR85b}
{\sc J.~C. Cox, J.~E. Ingersoll, Jr., and S.~A. Ross}, {\em A theory of the
  term structure of interest rates}, Econometrica, 53 (1985), pp.~385--407.

\bibitem{culina_stochastic_2011}
{\sc J.~Culina, S.~Kravtsov, and A.~H. Monahan}, {\em Stochastic
  {{Parameterization Schemes}} for {{Use}} in {{Realistic Climate Models}}},
  Journal of the Atmospheric Sciences, 68 (2011), pp.~284--299.

\bibitem{DubinkinaFrank07}
{\sc S.~Dubinkina and J.~Frank}, {\em Statistical mechanics of {A}rakawa's
  discretizations}, Journal of Computational Physics, 227 (2007), pp.~1286 --
  1305.

\bibitem{DubinkinaFrank10}
{\sc S.~Dubinkina and J.~Frank}, {\em Statistical relevance of vorticity
  conservation in the {H}amiltonian particle-mesh method}, Journal of
  Computational Physics, 229 (2010), pp.~2634 -- 2648.

\bibitem{E03}
{\sc W.~E}, {\em Analysis of the heterogeneous multiscale method for ordinary
  differential equations}, Comm. Math. Sci., 1 (2003), pp.~423--436.

\bibitem{EEtAl07}
{\sc W.~E, B.~Engquist, X.~Li, W.~Ren, and E.~Vanden-Eijnden}, {\em
  Heterogeneous multiscale methods: {A} review}, Comm. Comp. Phys., 2 (2007),
  pp.~367--450.

\bibitem{GearKevrekidis03}
{\sc C.~Gear and I.~Kevrekidis}, {\em Projective methods for differential
  equations}, SIAM J. Sci. Comp., 24 (2003), pp.~1091--1106.

\bibitem{Givonetal04}
{\sc D.~Givon, R.~Kupferman, and A.~Stuart}, {\em Extracting macroscopic
  dynamics: Model problems and algorithms}, Nonlinearity, 17 (2004),
  pp.~R55--127.

\bibitem{GottwaldEtAl17}
{\sc G.~Gottwald, D.~Crommelin, and C.~Franzke}, {\em {Ensemble-based
  Atmospheric Data Assimilation}}, in Nonlinear and Stochastic Climate
  Dynamics, C.~L.~E. Franzke and T.~J. O'Kane, eds., Cambridge University
  Press, Cambridge, 2017, pp.~209--240.

\bibitem{GottwaldMelbourne13}
{\sc G.~A. Gottwald and I.~Melbourne}, {\em {A Huygens principle for diffusion
  and anomalous diffusion in spatially extended systems}}, Proc. Natl. Acad.
  Sci. USA, 110 (2013), pp.~8411--8416.

\bibitem{GottwaldMelbourne13c}
{\sc G.~A. Gottwald and I.~Melbourne}, {\em Homogenization for deterministic
  maps and multiplicative noise}, Proceedings of the Royal Society A:
  Mathematical, Physical and Engineering Science, 469 (2013).

\bibitem{HaLuWa06}
{\sc E.~Hairer, C.~Lubich, and G.~Wanner}, {\em Geometric numerical
  integration: structure-preserving algorithms for ordinary differential
  equations}, vol.~31, Springer, 2006.

\bibitem{HorsthemkeBook}
{\sc W.~W. Horsthemke and R.~Lefever}, {\em Noise-induced transitions : theory
  and applications in physics, chemistry, and biology}, Springer series in
  synergetics, Springer-Verlag, Berlin, New York, 1984.

\bibitem{KellyMelbourne17}
{\sc D.~Kelly and I.~Melbourne}, {\em Deterministic homogenization for
  fast--slow systems with chaotic noise}, Journal of Functional Analysis, 272
  (2017), pp.~4063 -- 4102.

\bibitem{KevrekidisGearEtAl03}
{\sc I.~G. Kevrekidis, C.~W. Gear, J.~M. Hyman, G.~K. Panagiotis, O.~Runborg,
  and C.~Theodoropoulos}, {\em Equation-free, coarse-grained multiscale
  computation: {E}nabling microscopic simulators to perform system-level
  analysis}, Comm. Math. Sci., 1 (2003), pp.~715--762.

\bibitem{Khasminsky66}
{\sc R.~Z. Khasminsky}, {\em {On stochastic processes defined by differential
  equations with a small parameter}}, Theory of Probability and its
  Applications, 11 (1966), pp.~211--228.

\bibitem{Kurtz73}
{\sc T.~G. Kurtz}, {\em A limit theorem for perturbed operator semigroups with
  applications to random evolutions}, Journal of Functional Analysis, 12
  (1973), pp.~55--67.

\bibitem{LeimkuhlerReich}
{\sc B.~Leimkuhler and S.~Reich}, {\em {Simulating Hamiltonian Dynamics}},
  Cambridge University Press, Cambridge, 2005.

\bibitem{MTV02}
{\sc A.~Majda, I.~Timofeyev, and E.~Vanden-Eijnden}, {\em A priori tests of a
  stochastic mode reduction strategy}, Phys. D, 170 (2002), pp.~206--252.

\bibitem{MTV99}
{\sc A.~J. Majda, I.~Timofeyev, and E.~Vanden~Eijnden}, {\em Models for
  stochastic climate prediction}, Proceedings of the National Academy of
  Sciences, 96 (1999), pp.~14687--14691.

\bibitem{Majdaetal01}
{\sc A.~J. Majda, I.~Timofeyev, and E.~Vanden~Eijnden}, {\em A mathematical
  framework for stochastic climate models}, Communications on Pure and Applied
  Mathematics, 54 (2001), pp.~891--974.

\bibitem{Majda03}
{\sc A.~J. Majda, I.~Timofeyev, and E.~Vanden-Eijnden}, {\em {Systematic
  strategies for stochastic mode reduction in climate}}, Journal of the
  Atmospheric Sciences, 60 (2003), pp.~1705--1722.

\bibitem{MelbourneNicol05}
{\sc I.~Melbourne and M.~Nicol}, {\em Almost sure invariance principle for
  nonuniformly hyperbolic systems}, Commun. Math. Phys., 260 (2005),
  pp.~131--146.

\bibitem{MelbourneNicol08}
{\sc I.~Melbourne and M.~Nicol}, {\em Large deviations for nonuniformly
  hyperbolic systems}, Trans. Amer. Math. Soc., 360 (2008), pp.~6661--6676.

\bibitem{MelbourneNicol09}
{\sc I.~Melbourne and M.~Nicol}, {\em A vector-valued almost sure invariance
  principle for hyperbolic dynamical systems}, Annals of Probability, 37
  (2009), pp.~478--505.

\bibitem{MelbourneStuart11}
{\sc I.~Melbourne and A.~Stuart}, {\em A note on diffusion limits of chaotic
  skew-product flows}, Nonlinearity, 24 (2011), pp.~1361--1367.

\bibitem{monahan_stochastic_2011}
{\sc A.~H. Monahan and J.~Culina}, {\em Stochastic {{Averaging}} of {{Idealized
  Climate Models}}}, Journal of Climate, 24 (2011), pp.~3068--3088.

\bibitem{Papanicolaou76}
{\sc G.~C. Papanicolaou}, {\em Some probabilistic problems and methods in
  singular perturbations}, Rocky Mountain Journal of Mathematics, 6 (1976),
  pp.~653--674.

\bibitem{PavliotisStuart}
{\sc G.~A. Pavliotis and A.~M. Stuart}, {\em {Multiscale Methods: Averaging and
  Homogenization}}, Springer, New York, 2008.

\bibitem{TaoEtAl10}
{\sc M.~Tao, H.~Owhadi, and J.~E. Marsden}, {\em Nonintrusive and structure
  preserving multiscale integration of stiff {ODE}s, {SDE}s, and {H}amiltonian
  systems with hidden slow dynamics via flow averaging}, Multiscale Model.
  Simul., 8 (2010), pp.~1269--1324.

\end{thebibliography}

\end{document}